\documentclass[onecolumn]{autart}
\usepackage{}    
\usepackage[colorlinks,linkcolor=black,anchorcolor=red,citecolor=red]{hyperref}%

\usepackage{graphicx}           
\usepackage{dcolumn}            
\usepackage{bm}                 
\usepackage{bbm}
\usepackage{graphics}           
\usepackage{epsfig}             
\usepackage{times}              
\usepackage[fleqn]{amsmath}
\usepackage{amssymb}
\usepackage{extarrows}
\usepackage{amsfonts}
\usepackage{mathrsfs}
\usepackage{fancyhdr}
\usepackage{color}
\usepackage{subfigure}
\usepackage{enumerate}
\usepackage{epstopdf}
 \usepackage{subfiles}
   \graphicspath{{eps/}}
\usepackage [latin1]{inputenc} %
\allowdisplaybreaks[4] 

\newtheorem{theorem}{Theorem}  %
 
\newtheorem{corollary}[theorem]{Corollary}
\newtheorem{proposition}[theorem]{Proposition}

\newtheorem{remark}{Remark}

\begin{document}

\begin{frontmatter}

\title{Proof of Proposition 3.1 in   the paper titled  ``Backstepping control of a class of space-time-varying linear parabolic PDEs via time invariant kernel functions''}

\author{Qiaoling Chen$^{1}$}\ead{cql20202016@my.swjtu.edu.cn},
\author{Jun Zheng$^{1,2}$}\ead{zhengjun2014@aliyun.com},
\author{Guchuan Zhu$^{2}$}\ead{guchuan.zhu@polymtl.ca}

\address{$^{1}${School of Mathematics, Southwest Jiaotong University,
        Chengdu 611756, Sichuan, China}\\
        $^{2}$Department of Electrical Engineering, Polytechnique Montr\'{e}al, P.O. Box 6079, Station Centre-Ville, Montreal, QC, Canada H3T 1J4}

\begin{abstract}

We provide  a detailed proof  of  Proposition 3.1  in the paper  titled ``Backstepping control of a class of space-time-varying linear parabolic PDEs via time invariant kernel functions''.

 In the paper titled ``Backstepping control of a class of space-time-varying linear parabolic PDEs via time invariant kernel functions'', we addressed the problem of exponential stabilization and continuous dependence of solutions on initial data in different norms for a class of $1$-D linear parabolic PDEs with space-time-varying coefficients under backstepping boundary control. In order to stabilize the system without involving a Gevrey-like condition or the event-triggered scheme, a boundary feedback controller was designed via a time  invariant kernel function.  By using the approximative Lyapunov method, the exponential stability of the closed-loop system was established in the spatial $L^{p}$-norm and $W^{1,p}$-norm, respectively, whenever $p\in [1, +\infty]$. It was also shown that the solution to the considered system depends continuously on the spatial $L^{p}$-norm and $W^{1,p}$-norm, respectively, of the initial data.

\end{abstract}
\end{frontmatter}
  \section{{Brief review on the considered problem}}\label{Problem-formulation}

Before presenting the proof of Proposition 3.1 in the paper titled ``Backstepping control of a class of space-time-varying linear parabolic PDEs via time invariant kernel functions'', i.e., \cite{Chen:2023}, we would like to restate the considered problem.  Specifically, given certain initial data $w_0$, we study the problem of stabilization and continuous dependence  of solutions on  initial data  in  different norms for a class of  linear parabolic {PDEs} with   space and time dependent {coefficients}:
\begin{subequations}\label{original system}
\begin{align}
	 w_{t}(x,t)=&w_{xx}(x,t)+ c(x,t)w(x,t)+\int_0^xw(y,t)f(x,y)\text{d}y, (x,t)\in (0,1)\times {\mathbb{R}_{>0}},\label{1a}\\
	 w_x(0,t)=&0, t\in {\mathbb{R}_{>0}},\label{1b}\\
	 w_x(1,t)=&U(t),t\in {\mathbb{R}_{>0}}, \label{1c}\\
	 w(x,0)=&w_0(x),x\in(0,1),
\end{align}
\end{subequations}
where  $c: [0,1]\times {\mathbb{R}_{\geq0}} \rightarrow \mathbb{R} $ and $f:[0,1]\times[0,1]\rightarrow \mathbb{R} $     are given functions,   and  $U$ is the control input to be determined to stabilize the system.

We assume that
{
 $p\in  [1,+\infty)$ is either an arbitrary constant, or $p=+\infty$, {$\theta\in\left(0,\frac{1}{2}\right]$} is a constant,
 $f \in C([0,1]\times[0,1])$,
 }
{and the reaction coefficient} $c \in \mathcal{H}^{\theta,\frac{\theta}{2}}\left({\overline{Q}_{\infty}}\right) $   is given by
\begin{align*}
 c(x,t):=&c_1(x)+c_2(t)+c_3(x,t),
 \end{align*}
with some functions $ c_1\in C([0,1]),c_2\in C (\mathbb{R}_{\geq 0})$, and
$
 c_3 \in    {W^{2,1}_2(Q_T)}
$ for all $T\in \mathbb{R}_{>0}.$
 Assume further that
 \begin{align*}
 \sup_{(x,t)\in   {\overline{Q}_{\infty}}} c(x,t)<+\infty ,
\end{align*}
 and
     \begin{align*}
|c_3(x,t)-c_3(y,t)|\leq &L(t)   \omega(|x-y|) ,
\end{align*}
  holds true for all $x,y\in [0,1]$ and all $t\in \mathbb{R}_{\geq 0}$,  where   $\omega\in C([0,1])$ and $L\in C (\mathbb{R}_{\geq 0})$  are nonnegative functions satisfying
 \begin{align*}
\lim_{t\rightarrow+\infty} L(t)=0~~\text{and}~~
\omega(s)\leq 1,\forall s\in[0,1].
\end{align*}
Let  {$\lambda_0 >0$}   be an arbitrary   constant  and satisfy
$
\lambda_0>\sup_{(x,t)\in  \overline{Q}_{\infty} } c(x,t).
$
Let $D:=\{(x,y)\in \mathbb{R}^2|~0\leq y \leq x\leq 1\}$ and
\begin{align*}
 {\mu(x,y):=\lambda_0-c_1(x)+c_1(y)},\forall (x,y)\in D.
\end{align*}
Consider the following equation defined over $D$:
\begin{subequations}\label{kernel}
\begin{align}
 		  		  k_{xx}(x,y)-k_{yy}(x,y)
 	 = &  {\mu(x,y)}k(x,y)+f(x,y)
  +\int_y^xk(x,z)f(z,y)\text{d}z,  \label{kb}\\
   2\frac{\text{d}}{\text{d}x}\left(k(x,x)\right)=&{\lambda_0},\label{ka} \\
 		  k_{y}(x,0)=&0,\label{kc}\\
   k(0,0)=&0,\label{kd}
 	\end{align}
\end{subequations}
where   $\frac{\text{d}}{\text{d}x}\left(k(x,x)\right):=k_x(x,y)|_{y=x}+k_y(x,y)|_{y=x}$.

Then, by means of such a kernel function $k$, we define the boundary feedback control law as
\begin{align}
 U(t):=-k(1,1)w(1,t)-\int_{0}^{1}k_x(1,y)w(y,t)\text{d}y.\label{control law}
\end{align}
The existence, regularity and maximum estimate of a solution $k$ to the   equation~\eqref{kernel} is guaranteed by the following proposition, which is  Proposition~3.1 of \cite{Chen:2023}.  The proof will be given in Section~2.

 \begin{proposition}[{\cite[Proposition~3.1]{Chen:2023}}]\label{proposition kernel}
 The equation~\eqref{kernel}  admits a unique solution $k$  that belongs to $C^{2}\left(D\right)$ and satisfies
\begin{align}\label{maxk}
  |k(x,y)|\leq  Me^{2M},\forall (x,y)\in D,
   \end{align}
   where $M:=   \frac{ \overline{f}+{\overline{\Lambda}}} {2}$ with
   $
 {\overline{f}}:=    \max_{(x,y)\in[0,1]\times[0,1]}|f(x,y)|,$ and $
{{\overline{\Lambda}}}:=
 {\max\left\{\lambda_0,\max_{(x,y)\in\mathcal{D}}|\mu(x,y)| \right\}}.
$

\end{proposition}


Now, we define  the set
\begin{align*}
W_0:=\bigg\{&w|~w\in \mathcal{H}^{2+\theta}([0,1]),w_{x}(0)=0, w_{x}(1)=-k(1,1)w(1)-\int_{0}^{1}k_x(1,y)w(y)\text{d}y \bigg\},
\end{align*}
for the initial data of system \eqref{original system}, and denote
$
 \lambda(x,t):=\lambda_0-c\left(x,t\right),
$
 and
$
    {\underline{\lambda}}:=\inf_{(x,t)\in {\overline{Q}_{\infty}}}  \lambda(x,t).
$

The main result obtained in \cite{Chen:2023} is the following theorem, whose proof is given in \cite{Chen:2023}.
\begin{theorem}\label{theorem}
Given initial data {$w_0\in W_0$} and considering system \eqref{original system} under the {boundary} feedback control law \eqref{control law},  the following statements hold true:
 \begin{enumerate}[(i)]
\item   The system \eqref{original system}    admits a unique classical solution $w$, which belongs to
$C^{2,1} (\overline{Q}_T) $ for any $T\in \mathbb{R}_{>0}$.
 \item The system \eqref{original system}   is exponentially {stable  in the spatial} $L^p$-norm, having the  {estimate} for all {$t\in {\mathbb{R}_{>0}}$:}
\begin{align}\label{(ii)-p1}
\|w[t]\|_{L^p(0,1)}\leq C_{0}e^{-{\sigma_0} t}\|w_0\|_{L^p(0,1)},
\end{align}
 {where $\sigma_0$ and $C_{0}$ are positive constants depending only on $\underline{\lambda},k$ and $\max_{t\in \mathbb{R}_{\geq 0}}L(t)$ whenever $p\in[1,+\infty]$.}

 In addition, the solution to system \eqref{original system}  is continuously dependent on    the {spatial} $L^{p}$-norm   of the initial data,   having the {estimate} {for all $T\in \mathbb{R}_{>0}$:}
\begin{align*}
 &\|w_1-w_2\|_{C([0,T];L^p(0,1))}
 \leq C_0\| w_{01}-w_{02}\|_{L^p(0,1)},
 \end{align*}
  where $w_i  $ denotes the  solution to system \eqref{original system} corresponding to the  initial data $w_{0i}\in W_0  (i=1,2)$.
  \item Assume further that $c_1\in C^1 ([0,1]) $ and
      \begin{align*}
      \omega(s)=s^{\gamma_0},\forall s\in[0,1],
      \end{align*}
      with some constant $\gamma_0\geq 1$. Then, system \eqref{original system}   is exponentially {stable  in the spatial}   $W^{1,p}$-norm, having the estimate for all {$t\in {\mathbb{R}_{>0}}$:}
   \begin{align*}
\|w[t]\|_{W^{1,p}(0,1)}\leq \mathcal{C}_{0}e^{-\frac{\sigma_0}{2} t}\|w_0\|_{W^{1,p}(0,1)},
\end{align*}
where  {$\sigma_0$ is the same as in \eqref{(ii)-p1}, and $\mathcal{C}_{0}$
  is a positive constant  depending only on  $p,c_x,\underline{\lambda},k$ and $\max_{t\in \mathbb{R}_{\geq 0}}L(t)$ for $p\in[1,+\infty)$,  and only on $c_x,\underline{\lambda}$ and $k$ for $p=+\infty$, respectively.}

 In addition,  the solution to system \eqref{original system}  is continuously dependent on    the {spatial} $W^{1,p}$-norm   of the initial data,   having the estimate {for all $T\in \mathbb{R}_{>0}$:}
 \begin{align*}
  \|w_1-w_2\|_{C([0,T];W^{1,p}(0,1))}
 \leq&  \mathcal{C}_{0} \| w_{01}-w_{02}\|_{W^{1,p}(0,1)},
 \end{align*}
    where $w_i  $ denotes the  solution to system \eqref{original system} corresponding to the  initial data $w_{0i}\in W_0  (i=1,2)$.
  \end{enumerate}
\end{theorem}

{As a special case of Theorem~\ref{theorem},  the following result  indicates that, for  certain parabolic PDEs with time-varying coefficients, a boundary feedback control  can be  {designed} via  modified Bessel functions, which are often used to design boundary controls for stabilizing parabolic PDEs with constant coefficients.} The proof is also given in \cite{Chen:2023}.
\begin{corollary}\label{corollary}
 {{Let $c(x,t)\equiv c_2(t)$ and $f(x,y)\equiv  0$.} Under the following boundary control law\footnote{{Note that $
\lim_{s\rightarrow 0}\frac{\mathcal{I}_1(s)}{s} =\frac{1}{2}$ and $ \lim_{s\rightarrow 0}\frac{\mathcal{I}_2(s)}{s^2} =\frac{1}{8}$; {thus \eqref{special control law} is well-defined.}}}}
  \begin{align}
 U(t):=&-\frac{\lambda_0}{2}w(1,t)
  -\lambda_0 \int_{0}^{1} \frac{\mathcal{I}_1\left(\sqrt{\lambda_0(1-y^2)}\right)}{\sqrt{\lambda_0(1-y^2)}} w(y,t)\text{d}y -\lambda_0 \int_{0}^{1} \frac{\mathcal{I}_2\left(\sqrt{\lambda_0(1-y^2)}\right)}{  {1-y^2 } } w(y,t)\text{d}y,\label{special control law}
\end{align}
{system \eqref{original system}   is exponentially {stable  in the spatial} $L^p$-norm (and $W^{1,p}$-norm under the   further assumption that $c_2\in C^1(\mathbb{R}_{\geq 0})$), {where $\mathcal{I}_1$ and $\mathcal{I}_2$ are modified Bessel functions (see, e.g., \cite{Smyshlyaev:2004})  of order one and order two, respectively.}}
\end{corollary}



\section{Proof of  Proposition~\ref{proposition kernel}} \label{App.B}

In this section we provide a detailed proof for Proposition~\ref{proposition kernel}, i.e., {\cite[Proposition~3.1]{Chen:2023}}.




We prove  by using the method of   successive approximations as in, e.g., \cite{Colton:1977,krstic2008boundary,liu2003boundary,meurer2009tracking,si2018boundary,Smyshlyaev:2004,smyshlyaev2005control}. Indeed, {let}
$
\xi:=x+y$ and $\eta:=x-y.
$
{It follows that} $\eta\in[0,1]$ and $\xi\in[\eta,2-\eta]$. {Let}
\begin{align*}
G(\xi,\eta):=k(x,y)=k\left(\frac{\xi+\eta}{2},\frac{\xi-\eta}{2}\right).
\end{align*}
{Then} \eqref{kernel} is changed into
\begin{subequations}\label{Gset}
\begin{align}
4G_{\xi\eta}(\xi,\eta)
=& { \widetilde{\mu}(\xi,\eta)}G(\xi,\eta)+f \left(\frac{\xi+\eta}{2},\frac{\xi-\eta}{2}\right)
+\int_{\frac{\xi-\eta}{2}}^{\frac{\xi+\eta}{2}}
 \widetilde{\mathcal{G}}(z,\xi,\eta)\text{d}z, \label{Gb}\\
 G(\xi,0)=&{\frac{\lambda_0}{4}}\xi, \label{Ga}\\
G_\xi(\xi,\xi)=&G_\eta(\xi,\xi), \label{Gc}\\
G(0,0)=&0, \label{Gd}
\end{align}
\end{subequations}
where
\begin{align*}
 &{\widetilde{\mu}(\xi,\eta):=\mu(x,y)=\lambda_0-c_1\left(\frac{\xi+\eta}{2}\right) +c_1\left(\frac{\xi-\eta}{2}\right)},\\ &\widetilde{\mathcal{G}}(z,\xi,\eta):=f\left(z,\frac{\xi-\eta}{2}\right)G\left(\frac{\xi+\eta}{2}+z,\frac{\xi+\eta}{2}-z\right).
\end{align*}
Integrating \eqref{Gb} w.r.t. $\eta$ from $0$ to $\eta$ and differentiating \eqref{Ga} w.r.t. $\xi$ give
\begin{align}\label{Ge}
G_\xi(\xi,\eta)=&{\frac{\lambda_0}{4}}+{{\frac{{ 1}}{4}}\int_{0}^{\eta}\widetilde{\mu}(\xi,s)G(\xi,s)\text{d}s}+\frac{1}{4}\int_{0}^{\eta}f\left(\frac{\xi+s}{2},\frac{\xi-s}{2}\right)\text{d}s+\frac{1}{4}\int_{0}^{\eta}\int_{\xi}^{\xi+\eta-s}\widehat{G}(\tau,s,\xi)\text{d}\tau\text{d}s,
\end{align}
where $\widehat{G}(\tau,s,\xi):=f\left(\frac{\tau-s}{2},\xi-\frac{\tau+s}{2}\right)G(\tau,s)$.

Integrating  {\eqref{Ge}} w.r.t. $\xi$ from $\eta$ to $\xi$ gives
\begin{align}\label{equ.411-3}
G(\xi,\eta)=&G(\eta,\eta)+{\frac{\lambda_0}{4}}(\xi-\eta)+{{\frac{{ 1}}{4}}\int_{\eta}^{\xi}\int_{0}^{\eta}\widetilde{\mu}(\tau,s)G(\tau,s)\text{d}s\text{d}\tau}+\frac{1}{4}\int_{\eta}^{\xi}\int_{0}^{\eta}f\left(\frac{\tau+s}{2},\frac{\tau-s}{2}\right)\text{d}s\text{d}\tau\notag\\
&+\frac{1}{4}\int_{\eta}^{\xi}\int_{0}^{\eta}\int_{z}^{z+\eta-s}\widehat{G}(\tau,s,z)\text{d}\tau\text{d}s\text{d}z.
\end{align}
By \eqref{Gc}, we have
\begin{align}
\frac{\text{d}}{\text{d}\xi}\left(G(\xi,\xi)\right)=2G_\xi(\xi,\xi).\label{equ.411-1}
\end{align}
  Substituting $\eta=\xi$ into \eqref{Ge}, {then integrating it over $[0,\xi]$ and applying \eqref{Gd} and \eqref{equ.411-1}, we obtain}
\begin{align*}
G(\xi,\xi)=&{\frac{\lambda_0}{2}}\xi+{{\frac{{ 1}}{2}}\int_{0}^{\xi}\int_{0}^{\tau}\widetilde{\mu}(\tau,s)G(\tau,s)\text{d}s\text{d}\tau}+\frac{1}{2}\int_{0}^{\xi}\int_{0}^{\tau}f\left(\frac{\tau+s}{2},\frac{\tau-s}{2}\right)\text{d}s\text{d}\tau\notag\\
&+\frac{1}{2}\int_{0}^{\xi}\int_{0}^{z}\int_{z}^{2z-s}\widehat{G}(\tau,s,z)\text{d}\tau\text{d}s\text{d}z,
\end{align*}
{which gives
\begin{align}\label{equ.411-2}
G(\eta,\eta)=&{\frac{\lambda_0}{2}}\eta+{{\frac{{ 1}}{2}}\int_{0}^{\eta}\int_{0}^{\tau}\widetilde{\mu}(\tau,s)G(\tau,s)\text{d}s\text{d}\tau}+\frac{1}{2}\int_{0}^{\eta}\int_{0}^{\tau}f\left(\frac{\tau+s}{2},\frac{\tau-s}{2}\right)\text{d}s\text{d}\tau\notag\\
&+\frac{1}{2}\int_{0}^{\eta}\int_{0}^{z}\int_{z}^{2z-s}\widehat{G}(\tau,s,z)\text{d}\tau\text{d}s\text{d}z.
\end{align}}
{We deduce by \eqref{equ.411-3} and \eqref{equ.411-2} that}
\begin{align}\label{otherkernel}
G(\xi,\eta)=&{\frac{\lambda_0}{4}}\left(\xi+\eta\right)
+\frac{1}{4}\int_{\eta}^{\xi}\int_{0}^{\eta}f\left(\frac{\tau+s}{2},\frac{\tau-s}{2}\right)\text{d}s\text{d}\tau+\frac{1}{2}\int_{0}^{\eta}\int_{0}^{\tau}f\left(\frac{\tau+s}{2},\frac{\tau-s}{2}\right)\text{d}s\text{d}\tau\notag\\
&+{{\frac{{ 1}}{4}}\int_{\eta}^{\xi}\int_{0}^{\eta}\widetilde{\mu}(\tau,s)G(\tau,s)\text{d}s\text{d}\tau}+{{\frac{{ 1}}{2}}\int_{0}^{\eta}\int_{0}^{\tau}\widetilde{\mu}(\tau,s)G(\tau,s)\text{d}s\text{d}\tau}
+\frac{1}{4}\int_{\eta}^{\xi}\int_{0}^{\eta}\int_{z}^{z+\eta-s}\widehat{G}(\tau,s,z)\text{d}\tau\text{d}s\text{d}z\notag\\
&+\frac{1}{2}\int_{0}^{\eta}\int_{0}^{z}\int_{z}^{2z-s}\widehat{G}(\tau,s,z)\text{d}\tau\text{d}s\text{d}z.
\end{align}
We rewrite the integral equation \eqref{otherkernel}   as
\begin{align}\label{twoparts}
G(\xi,\eta)=G_0(\xi,\eta)+\Phi_G(\xi,\eta),
\end{align}
where  we defined
\begin{align*}
G_0(\xi,\eta):=&{\frac{\lambda_0}{4}}\left(\xi+\eta\right)
+\frac{1}{4}\int_{\eta}^{\xi}\int_{0}^{\eta}f\left(\frac{\tau+s}{2},\frac{\tau-s}{2}\right)\text{d}s\text{d}\tau
+\frac{1}{2}\int_{0}^{\eta}\int_{0}^{\tau}f\left(\frac{\tau+s}{2},\frac{\tau-s}{2}\right)\text{d}s\text{d}\tau,
\end{align*}
{and}  $\Phi_G(\xi,\eta)$   via
\begin{align}
{\Phi_H(\xi,\eta)} :=&{{\frac{{ 1}}{4}}\int_{\eta}^{\xi}\int_{0}^{\eta}\widetilde{\mu}(\tau,s)H(\tau,s)\text{d}s\text{d}\tau}
+{{\frac{{ 1}}{2}}\int_{0}^{\eta}\int_{0}^{\tau}\widetilde{\mu}(\tau,s)H(\tau,s)\text{d}s\text{d}\tau}
+\frac{1}{4}\int_{\eta}^{\xi}\int_{0}^{\eta}\int_{z}^{z+\eta-s}\widehat{H}(\tau,s,z)\text{d}\tau\text{d}s\text{d}z\notag\\
&+\frac{1}{2}\int_{0}^{\eta}\int_{0}^{z}\int_{z}^{2z-s}\widehat{H}(\tau,s,z)\text{d}\tau\text{d}s\text{d}z,\notag\\
\widehat{H}(\xi,\eta,z):=&f\left(\frac{\xi-\eta}{2},z-\frac{\xi+\eta}{2}\right)H(\xi,\eta),\label{hatH}
\end{align}
for  any  $\eta\in[0,1]$, $\xi\in[\eta,2-\eta]$, {$z\in[\eta,2-\eta]$}, and any function $H$.

Define the sequence  {$\{G_n(\xi,\eta)\}$} {for} $n\in \mathbb{N}_0$ via
\begin{align}\label{sequence}
G_{n+1}(\xi,\eta):=G_0(\xi,\eta)+{\Phi_{G_n}}(\xi,\eta),
\end{align}
and {let}
$\Delta G_{n}(\xi,\eta):=G_{n+1}(\xi,\eta)-G_{n}(\xi,\eta)$. It holds that
\begin{align*}
\Delta G_{n+1}(\xi,\eta)={\Phi_{\Delta G_n}}(\xi,\eta)
\end{align*}
and
\begin{align}\label{series sequence}
G_{n+1}{(\xi,\eta)}=G_0{(\xi,\eta)}+\sum_{j=0}^{n}\Delta G_{j}(\xi,\eta).
\end{align}
From \eqref{sequence}, we {see that if $\{G_n(\xi,\eta)\}$ converges uniformly w.r.t. $(\xi,\eta)$} when $n\rightarrow\infty$, then  $G(\xi,\eta):=\lim_{n\to\infty}G_n(\xi,\eta)$ is the solution {to} the integral equation {\eqref{twoparts}}, {and therefore the existence of a kernel function $k$ {given by} \eqref{kernel} is guaranteed. In addition, in view of} \eqref{series sequence}, {the convergence of  $\{G_n(\xi,\eta)\}$ is equivalent to the convergence of the series $\sum_{n=0}^{\infty}\Delta G_n(\xi,\eta)$. {Therefore,} it suffices} to show that the series $\sum_{n=0}^{\infty}\Delta G_n(\xi,\eta)$ is uniformly convergent {w.r.t. {$(\xi,\eta)$}. Furthermore, recalling  the Weierstrass M-test, it suffices to show that  the following inequality}
\begin{align}\label{induction}
\left|\Delta G_n(\xi,\eta)\right| \leq \frac{M^{n+2}}{(n+1)!}\left(\xi+\eta\right)^{n+1}
\end{align}
{holds true for all $n\in\mathbb{N}_0$,} where
\begin{align*}
M:=&  \frac{ \overline{f}+{\overline{\Lambda}}} {2},
 {\overline{f}}:=   \max_{(x,y)\in[0,1]\times [0,1]}|f(x,y)|,
{{\overline{\Lambda}}}:=
 {\max\left\{\lambda_0,\max_{(x,y)\in\mathcal{D}}|\lambda_0-c_1(x)+c_1(y)| \right\}}.
\end{align*}
 {Now we prove  \eqref{induction}  by induction. First of all, since}
{\begin{align}\label{G0}
\left| G_0(\xi,\eta)\right|&\leq{\frac{\lambda_0}{4}}\left(\xi+\eta\right)+\frac{{\overline{f}}}{4}(\xi\eta-\eta^2)+\frac{{\overline{f}}}{4}\eta^2
 \leq M,
\end{align}}
 {it follows that}
\begin{align}
\left|\Delta G_{0}(\xi,\eta)\right|=&\left|\Phi _{ G_{0}}(\xi,\eta) \right|\notag\\
\leq&\frac{{{{\overline{\Lambda}}}}}{4}\int_{\eta}^{\xi}\int_{0}^{\eta}M\text{d}s\text{d}\tau\notag +{\frac{{{{\overline{\Lambda}}}}}{2}}\int_{0}^{\eta}\int_{0}^{\tau}M\text{d}s\text{d}\tau
+\frac{\overline{f}}{4}\int_{\eta}^{\xi}\int_{0}^{\eta}\int_{z}^{z+\eta-s}M\text{d}\tau\text{d}s\text{d}z
+\frac{\overline{f}}{2}\int_{0}^{\eta}\int_{0}^{z}\int_{z}^{2z-s}M\text{d}\tau\text{d}s\text{d}z\notag\\
=&{M\left(\frac{{{{\overline{\Lambda}}}}}{4}\xi\eta+\frac{{\overline{f}}}{8}\xi\eta^2-\frac{{\overline{f}}}{24}\eta^3\right)}\notag\\
\leq&{M^2\left(\xi+\eta\right),}
\end{align}
which shows that \eqref{induction} holds true for $n=0$.

{Supposing that the inequality \eqref{induction}  holds true for a general $n\in \mathbb{N}$, we need to show that \eqref{induction} also holds true for  $n+1$. Indeed,} we get
\begin{align*} 
\left|\Delta G_{n+1}(\xi,\eta)\right|=&\left|{\Phi_{\Delta G_n}}(\xi,\eta)  \right|
=  |I_1+I_2+I_3+I_4|,
\end{align*}
where
\begin{align*}
I_1:=&{\frac{{ 1}}{4}}\int_{\eta}^{\xi}\int_{0}^{\eta}{\widetilde{\mu}(\tau,s)}\Delta G_{n}(\tau,s)\text{d}s\text{d}\tau,\notag\\
I_2:=&{\frac{{ 1}}{2}}\int_{0}^{\eta}\int_{0}^{\tau}{\widetilde{\mu}(\tau,s)}\Delta G_{n}(\tau,s)\text{d}s\text{d}\tau,\notag\\
I_3:=&\frac{1}{4}\int_{\eta}^{\xi}\int_{0}^{\eta}\int_{z}^{z+\eta-s}\widehat{\Delta G_{n}}(\tau,s,z)\text{d}\tau\text{d}s\text{d}z,\\
I_4:=&\frac{1}{2}\int_{0}^{\eta}\int_{0}^{z}\int_{z}^{2z-s}\widehat{\Delta G_{n}}(\tau,s,z)\text{d}\tau\text{d}s\text{d}z,
\end{align*}
 {with} $\widehat{\Delta G_{n}} $ defined via \eqref{hatH}.

{Since  \eqref{induction} holds true for $n$, it {follows} that}
{\begin{align*}
|I_1|\leq&\frac{{{\overline{\Lambda}}}M^{n+2}}{4(n+1)!}\int_{\eta}^{\xi}\int_{0}^{\eta}(\tau+s)^{n+1}\text{d}s\text{d}\tau
\leq\frac{{{\overline{\Lambda}}}M^{n+2}}{4(n+1)!}\frac{1}{(n+2)(n+3)}\left((\xi+\eta)^{n+3}-(2\eta)^{n+3}\right),\\
 |I_2|\leq&\frac{{{\overline{\Lambda}}}M^{n+2}}{2(n+1)!}\int_{0}^{\eta}\int_{0}^{\tau}(\tau+s)^{n+1}\text{d}s\text{d}\tau
\leq \frac{{{\overline{\Lambda}}}M^{n+2}}{4(n+1)!}\frac{1}{(n+2)(n+3)}(2\eta)^{n+3},\\
|I_3|\leq&\frac{\overline{f} M^{n+2}}{4(n+1)!}\int_{\eta}^{\xi}\int_{0}^{\eta}\int_{z}^{z+\eta-s}(\tau+s)^{n+1}\text{d}\tau\text{d}s\text{d}z
\leq \frac{\overline{f} M^{n+2}}{4(n+1)!}\frac{1}{(n+2)(n+3)}\left((\xi+\eta)^{n+3}-(2\eta)^{n+3}\right),\\
|I_4|\leq&\frac{\overline{f} M^{n+2}}{2(n+1)!}\int_{0}^{\eta}\int_{0}^{z}\int_{z}^{2z-s}(\tau+s)^{n+1}\text{d}\tau\text{d}s\text{d}z
\leq \frac{\overline{f} M^{n+2}}{2(n+1)!}\frac{1}{(n+2)(n+3)}(2\eta)^{n+3}.
\end{align*}}
Then we get
\begin{align*}
	 |\Delta G_{n+1}(\xi,\eta)|
\leq &	\left(|I_1|+|I_2|\right)+\left(|I_3|+|I_4|\right)\\
\leq &\frac{{{{\overline{\Lambda}}}}}{4(n+3)}\frac{M^{n+2}}{(n+2)!}(\xi+\eta)^{n+3}
{+}\frac{{\overline{f}}}{2(n+3)}\frac{M^{n+2}}{(n+2)!}(\xi+\eta)^{n+3}\\
\leq & \left(\frac{{{{\overline{\Lambda}}}}}{2(n+3)}+\frac{{\overline{f}}}{n+3}\right)\frac{M^{n+2}}{(n+2)!}(\xi+\eta)^{n+2}\\
\leq &\frac{M^{n+3}}{(n+2)!}\left(\xi+\eta\right)^{n+2},
	\end{align*}
which implies that \eqref{induction} holds for $n+1$. Therefore, \eqref{induction} holds {true} for $n\in \mathbb{N}_0$.

 {We conclude} that the series
$
\sum_{n=0}^\infty \Delta G_n(\xi,\eta)
$
converges absolutely and uniformly in $\eta\in[0,1]$ and $\xi\in[\eta,2-\eta]$, and that
\begin{align}\label{G-lim}
G(\xi,\eta):=\lim_{n\to\infty}G_n(\xi,\eta)
\end{align}
  is a  solution to {\eqref{twoparts}. Moreover,  $G(\xi,\eta)$ is continuous in {$(\xi,\eta)$}}. Furthermore, in view of \eqref{otherkernel} with $f\in C  \left([0,1]\times[0,1]\right)$ and $c_1\in C \left([0,1] \right)$,   we deduce that $G$ is $C^2$-continuous.  By using the change of variables, we deduce that the problem  \eqref{kernel} admits a  solution   $k\in C^2\left(D\right)$.

Regarding the  uniqueness of solution to \eqref{kernel}, which is equivalent to the uniqueness of  solution to the integral equation \eqref{otherkernel},  it suffices to show that any twice continuously differentiable solution  of  \eqref{otherkernel} can be approximated by the sequence $\{G_n(\xi,\eta)\}$ constructed via \eqref{sequence}. More precisely, for  any twice continuously differentiable solution   $\overline{G}(\xi,\eta)$ of \eqref{otherkernel}, we intend to show that the estimate
\begin{align} \label{unique}
|G_n(\xi,\eta)-\overline{G}(\xi,\eta)|\leq \frac{\delta M^{n}}{n!}\left(\xi+\eta\right)^{n}
\end{align}
holds true for all $n\in\mathbb{N}_0$, where
\begin{align*}
 \delta:= \max_{\eta\in[0,1],\xi\in[\eta,2-\eta]} |\Phi_{\overline{G}}(\xi,\eta) |.
\end{align*}
This is because  \eqref{unique} implies the uniform  convergence of $\{G_n(\xi,\eta)\}$ to $ \overline{G}(\xi,\eta)$, while   the uniqueness of $ G(\xi,\eta):= \lim_{n\to\infty}G_n(\xi,\eta)$  leads to $\overline{G}(\xi,\eta)\equiv G(\xi,\eta) $. {Thus,} the solution {to} \eqref{otherkernel} is
unique, and consequently, the solution to  {\eqref{kernel}} is also unique.

Now we prove \eqref{unique} by induction. First, since $\overline{G}(\xi,\eta)$ is a solution {to} \eqref{otherkernel}, using the definition of $\Phi_{\overline{G}}$, we obtain
\begin{align*}
|G_0(\xi,\eta)-\overline{G}(\xi,\eta)|
 =|G_0(\xi,\eta)-\left(G_0(\xi,\eta)+ \Phi_{\overline{G}}(\xi,\eta)\right)|
 =|\Phi_{\overline{G}}(\xi,\eta)|
 \leq\delta,
\end{align*}
which shows that \eqref{unique} holds true for $n=0$.

Suppose that \eqref{unique} holds true for a general $n\in \mathbb{N}$. Then, by  this assumption and the linearity of $\Phi_{H}(\xi,\eta)$  w.r.t. $H$, we have
\begin{align*}
|G_{n+1}(\xi,\eta)-\overline{G}(\xi,\eta)|=\left|     \Phi_{G_n}(\xi,\eta)  -\Phi_{\overline{G}}(\xi,\eta)\right |
=\left|     \Phi_{G_n -\overline{G}}(\xi,\eta)\right |
\leq I_5+I_6+I_7+I_8,
\end{align*}
where
\begin{align*}
I_5:=&\delta\frac{{{\overline{\Lambda}}}M^{n}}{4n!}\int_{\eta}^{\xi}\int_{0}^{\eta}(\tau+s)^{n}\text{d}s\text{d}\tau,\notag\\
I_6:=&\delta\frac{{{\overline{\Lambda}}}M^{n}}{2n!}\int_{0}^{\eta}\int_{0}^{\tau}(\tau+s)^{n}\text{d}s\text{d}\tau,\notag\\
I_7:=&\delta\frac{\overline{f} M^{n}}{4n!}\int_{\eta}^{\xi}\int_{0}^{\eta}\int_{z}^{z+\eta-s}(\tau+s)^{n}\text{d}\tau\text{d}s\text{d}z,\notag\\
I_8:=&\delta\frac{\overline{f} M^{n} }{2n!}\int_{0}^{\eta}\int_{0}^{z}\int_{z}^{2z-s}(\tau+s)^{n}\text{d}\tau\text{d}s\text{d}z.
\end{align*}
Analogous to the estimates of $I_1,I_2,I_3$ and $I_4$, it holds that
\begin{align*}
I_5 \leq&\delta\frac{{{\overline{\Lambda}}}M^{n}}{4n!}\frac{1}{(n+1)(n+2)}\left((\xi+\eta)^{n+2}-(2\eta)^{n+2}\right),\\
I_6 \leq&\delta\frac{{{\overline{\Lambda}}}M^{n}}{4n!}\frac{1}{(n+1)(n+2)}(2\eta)^{n+2},\\
I_7\leq&\delta\frac{\overline{f} M^{n}}{4n!}\frac{1}{(n+1)(n+2)}\left((\xi+\eta)^{n+2}-(2\eta)^{n+2}\right),\\
I_8 \leq&\delta\frac{\overline{f} M^{n} }{2n!}\frac{1}{(n+1)(n+2)}(2\eta)^{n+2}.
	\end{align*}
Then, we obtain
\begin{align*}
	 |G_{n+1}(\xi,\eta)-\overline{G}(\xi,\eta)|
	\leq  \delta \left(\frac{{{{\overline{\Lambda}}}}}{2(n+2)}+\frac{{\overline{f}}}{n+2}\right)\frac{M^{n}}{(n+1)!}(\xi+\eta)^{n+1}
\leq  {\frac{\delta M^{n+1}}{(n+1)!}\left(\xi+\eta\right)^{n+1}},
	\end{align*}
which implies that \eqref{unique} holds true for $n+1$. Thus,  \eqref{unique} holds true for all $n\in\mathbb{N}_0$. We conclude that the solution {to} \eqref{kernel} is unique.

Finally,  the estimate \eqref{maxk} follows immediately from  \eqref{G-lim}, \eqref{G0},  \eqref{series sequence}, \eqref{induction} and  the change of variables. The proof is {completed}.

\begin{remark}\label{Remark 2}
For the stabilization problem of linear parabolic PDEs with time  invariant coefficients, the ``uniqueness  of kernel functions'' was claimed in \cite{Smyshlyaev:2004} and proved by estimating
\begin{align}
 |\Delta G(\xi,\eta)|:= |G'(\xi,\eta)-G''(\xi,\eta)|\leq \frac{C_0}{n!}(\xi+\eta)^n\label{error}
\end{align}
with some constant $C_0>0$, where $G',G''$ are two kernel functions; see (39) and (40) of \cite{Smyshlyaev:2004}.

It is worth noting that a gap exists in the proof {given in} \cite{Smyshlyaev:2004}.
On one hand,  if $G' $ and $G''$  are given by (see (39) of \cite{Smyshlyaev:2004})
\begin{align}\label{47}
 G'(\xi,\eta)=\sum_{n=0}^{\infty}G_n(\xi,\eta), ~
 G''(\xi,\eta)=\sum_{n=0}^{\infty}G_n(\xi,\eta),
\end{align}
as indicated in  \cite{Smyshlyaev:2004}, the estimate  \eqref{error} can be obtained in the same way as in (38) of \cite{Smyshlyaev:2004}.
Therefore  the equality of  $G'=G''$ can be concluded by  letting $n\rightarrow \infty$. However, {since the limit of $\sum_{i=0}^{n}G_i(\xi,\eta)$ is unique, which is known in the limit theory, due to \eqref{47}}, the equality of $G'=G''$ holds true trivially. {Therefore,} the proof presented in \cite{Smyshlyaev:2004}
is essentially for the uniqueness of the limit of $\{G_n\}$ rather than the uniqueness of the solution.

On the other hand, if $G'$ and $G''$ are not given by  \eqref{47}, {without imposing more conditions}, the estimate  \eqref{error} cannot be obtained in the same way as in (38) of \cite{Smyshlyaev:2004}.

In order to  address the {aforementioned} {issues}, as presented in this paper, a feasible idea for proving the uniqueness of kernel functions  is to prove that all solutions to integral equations ({e.g., \eqref{otherkernel} or \eqref{twoparts})} can be approximated by the sequence  {$\{G_n\}$}. {That is to say, $G'$ and $G''$ in \cite{Smyshlyaev:2004}    have to be with a form of \eqref{47}}.    Then  the uniqueness of the solution {is determined} by the uniqueness of the limit of $\{G_n\}$. Note that,  by using this method, the estimate considered in this paper (see \eqref{unique}) is   different from  \eqref{error}. In particular, \eqref{unique} implies \eqref{error}.
\end{remark}

\begin{thebibliography}{}
\bibitem{Chen:2023}
Q. Chen, J. Zheng, and G. Zhu.
\newblock {\em Backstepping control of a class of
space-time-varying linear parabolic PDEs via time invariant kernel
functions}.
\newblock  2023.

%
%
%

%
%
%

\bibitem{Colton:1977}
D.~Colton.
\newblock The solution of initial-boundary value problems for parabolic
  equations by the method of integral operators.
\newblock {\em J. Diff. Equat.}, 26(2):181--190, 1977.

%
%
%
%
%
%
%

%
%
%
%
%
%
%
%
%
%

\bibitem{krstic2008boundary}
M.~Krstic and A.~Smyshlyaev.
\newblock {\em Boundary Control of {PDEs}: A Course on Backstepping Designs}.
\newblock SIAM, 2008.

%
%
%


\bibitem{liu2003boundary}
W.~Liu.
\newblock Boundary feedback stabilization of an unstable heat equation.
\newblock {\em SIAM J. Control Optim.}, 42(3):1033--1043, 2003.

\bibitem{meurer2009tracking}
T.~Meurer and A.~Kugi.
\newblock Tracking control for boundary controlled parabolic {PDEs} with
  varying parameters: Combining backstepping and differential flatness.
\newblock {\em Automatica}, 45(5):1182--1194, 2009.

%
%
%
%
%
\bibitem{si2018boundary}
Y.~Si, C.~Xie, and N.~Zhao.
\newblock Boundary control for a class of reaction-diffusion systems.
\newblock {\em Int. J. Autom. Comput.}, 15(1):94--102, 2018.

\bibitem{Smyshlyaev:2004}
A.~Smyshlyaev and M.~Krstic.
\newblock Closed-form boundary state feedbacks for a class of {1-D} partial
  integro-differential equations.
\newblock {\em IEEE Trans. Autom. Contr.}, 34(12):435--443, 2004.

%

\bibitem{smyshlyaev2005control}
A.~Smyshlyaev and M.~Krstic.
\newblock On control design for {PDEs} with space-dependent diffusivity or
  time-dependent reactivity.
\newblock {\em Automatica}, 41(9):1601--1608, 2005.


%
%
%
%

\end{thebibliography}
\end{document}